\documentclass[a4paper,reqno,10pt]{amsart}
\usepackage{amssymb,graphicx,latexsym}

\newcommand{\sym}{\mathcal{S}} 
\newcommand{\pt}{\sigma} 
\newcommand{\T}{\mathcal{T}} 
\newcommand{\R}{\mathcal{R}} 
\newcommand{\N}{{\mathbb N}}

\DeclareMathOperator{\proj}{\mathrm{proj}}
\DeclareMathOperator{\des}{des}
\DeclareMathOperator{\res}{Res}
\DeclareMathOperator{\ddif}{Ddif}
\newcommand{\pq}[1]{[\,#1\,]_{p,q}}
\newcommand{\p}[1]{[\,#1\,]_{p}}

\theoremstyle{plain} 
\newtheorem{theorem}{Theorem}
\newtheorem*{theorem*}{Theorem}
\newtheorem{corollary}[theorem]{Corollary}
\newtheorem*{corollary*}{Corollary}
\newtheorem{lemma}[theorem]{Lemma}
\newtheorem*{lemma*}{Lemma}
\newtheorem{proposition}[theorem]{Proposition}
\newtheorem*{proposition*}{Proposition}

\newtheorem*{thmu1}{Theorem \ref{u1}$^{\prime}$}
\newtheorem*{thmv1}{Theorem \ref{v1}$^{\prime}$}

\theoremstyle{definition} 

\newtheorem*{definition*}{Definition}
\newtheorem{example}[theorem]{Example} 
\newtheorem*{example*}{Example}

\theoremstyle{remark} 

\newtheorem*{remark*}{Remark}

\def\dd{\makebox[1.2ex]{\rule[.6ex]{0.8ex}{.15ex}}}

\def\ab{\ensuremath{12}}
\def\ba{\ensuremath{21}}

\def\abc{\ensuremath{123}}

\def\bac{\ensuremath{213}}

\def\axbc{\ensuremath{1{\dd}23}}
\def\axcb{\ensuremath{1{\dd}32}}
\def\bxac{\ensuremath{2{\dd}13}}
\def\bxca{\ensuremath{2{\dd}31}}
\def\cxab{\ensuremath{3{\dd}12}}
\def\cxba{\ensuremath{3{\dd}21}}

\def\abxc{\ensuremath{12{\dd}3}}
\def\acxb{\ensuremath{13{\dd}2}}
\def\baxc{\ensuremath{21{\dd}3}}
\def\bcxa{\ensuremath{23{\dd}1}}
\def\caxb{\ensuremath{31{\dd}2}}
\def\cbxa{\ensuremath{32{\dd}1}}

\def\axbxc{\ensuremath{1{\dd}2{\dd}3}}

\def\bxaxc{\ensuremath{2{\dd}1{\dd}3}}

\begin{document}

\title[Counting patterns of type $(1,2)$ or $(2,1)$]
{Counting occurrences of a pattern of type \\$(1,2)$ or $(2,1)$ in permutations}

\author{Anders Claesson}
\address{Matematik\\
  Chalmers tekniska h\"ogskola och G\"oteborgs universitet\\
  S-412 96 G\"oteborg, Sweden} 
\email{claesson@math.chalmers.se}

\author{Toufik Mansour} 
\address{
  LaBRI\\
  Univerit\'e Bordeaux I\\
  351 cours de la Lib\'eration, 33405 Talence Cedex, France }
\email{toufik@labri.fr}

\date{\today}

\begin{abstract}
  Babson and Steingr\'{\i}msson introduced generalized permutation
  patterns that allow the requirement that two adjacent letters in a
  pattern must be adjacent in the permutation. Claesson presented a
  complete solution for the number of permutations avoiding any single
  pattern of type $(1,2)$ or $(2,1)$.  For eight of these twelve
  patterns the answer is given by the Bell numbers. For the remaining
  four the answer is given by the Catalan numbers.
 
  With respect to being equidistributed there are three different
  classes of patterns of type $(1,2)$ or $(2,1)$. We present a
  recursion for the number of permutations containing exactly one
  occurrence of a pattern of the first or the second of the
  aforementioned classes, and we also find an ordinary generating
  function for these numbers. We prove these results both
  combinatorially and analytically. Finally, we give the distribution
  of any pattern of the third class in the form of a continued fraction,
  and we also give explicit formulas for the number of permutations
  containing exactly $r$ occurrences of a pattern of the third class
  when $r\in\{1,2,3\}$.
\end{abstract}

\maketitle\thispagestyle{empty}

\section{Introduction and preliminaries}

Let $[n]=\{1,2,\ldots,n\}$ and denote by $\sym_n$ the set of
permutations of $[n]$. We shall view permutations in $\sym_n$ as words
with $n$ distinct letters in $[n]$.

Classically, a pattern is a permutation
$\sigma\in\sym_k$, and an occurrence of $\sigma$ in a permutation
$\pi=a_1 a_2 \cdots a_n\in\sym_n$ is a subword of $\pi$ that is order
equivalent to $\sigma$. For example, an occurrence of $132$ is a
subword $a_i a_j a_k$ ($1\leq i < j < k\leq n$) of $\pi$ such that
$a_i < a_k < a_j$. We denote by $s_{\sigma}^r(n)$ the number of
permutations in $\sym_n$ that contain exactly $r$ occurrences of the
pattern $\sigma$.

In the last decade much attention has been paid to the problem of
finding the numbers $s_{\sigma}^r(n)$ for a fixed $r\geq 0$ and a
given pattern $\tau$ (see
\cite{AlFr00,At99,Bo97b,Bo97a,Bo98,ChWe99,ElNo01,Ma00,MaVa01,NoZe96,
  Ro99,SiSc85,St94,St96,We95}
). 
Most of the authors consider only the case $r=0$, thus studying
permutations {\it avoiding\/} a given pattern.  Only a few papers
consider the case $r>0$, usually restricting themselves to
patterns of length $3$.  Using two simple involutions (\emph{reverse}
and \emph{complement}) on $\sym_n$ it is immediate that with respect
to being equidistributed, the six patterns of length three fall into
the two classes $\{123,321\}$ and $\{132,213,231,312\}$. Noonan
\cite{No96} proved that $s_{123}^1(n)=\frac 3n\binom{2n}{n-3}$. A
general approach to the problem was suggested by Noonan and Zeilberger
\cite{NoZe96}; they gave another proof of Noonan's result, and
conjectured that
$$
s_{123}^2(n)=\frac{59n^2+117n+100}{2n(2n-1)(n+5)}\binom{2n}{n-4}
$$ 
and $s_{132}^1(n)=\binom{2n-3}{n-3}$. The latter conjecture was proved by
B\'ona in \cite{Bo98}. A conjecture of Noonan and Zeilberger states
that $s_{\sigma}^r(n)$ is $P$-recursive in $n$ for any $r$ and $\tau$. 
It was proved by B\'ona \cite{Bo97c} for $\sigma=132$. 

Mansour and Vainshtein \cite{MaVa01} suggested a new approach
to this problem in the case $\sigma=132$, which allows one to get an
explicit expression for $s_{132}^r(n)$ for any given $r$. More
precisely, they presented an algorithm that computes the generating
function $\sum_{n\geq0} s_{132}^r(n)x^n$ for any $r\geq 0$. To get
the result for a given $r$, the algorithm performs certain routine
checks for each element of the symmetric group $S_{2r}$. The algorithm
has been implemented in C, and yields explicit results for $1\leq
r\leq 6$.

In \cite{BaSt00} Babson and Steingr\'{\i}msson introduced generalized
permutation patterns that allow the requirement that two adjacent
letters in a pattern must be adjacent in the permutation.  The
motivation for Babson and Steingr\'{\i}msson in introducing these
patterns was the study of Mahonian permutation statistics. Two
examples of (generalized) patterns are $\axcb$ and $\acxb$. An
occurrence of $\axcb$ in a permutation $\pi = a_1 a_2 \cdots a_n$ is a
subword $a_i a_j a_{j+1}$ of $\pi$ such that $a_i<a_{j+1}<a_j$.
Similarly, an occurrence of $\acxb$ is a subword $a_i a_{i+1} a_{j}$
of $\pi$ such that $a_i<a_{j}<a_{i+1}$. More generally, if
$xyz\in\sym_3$ and $\pi = a_1 a_2 \cdots a_n\in\sym_n$, then we define
$$
(x\dd yz)\,\pi \,=\,
|\{ a_i a_j a_{j+1} : \proj(a_i a_j a_{j+1}) = xyz, 1\leq i < j < n \}|,
$$
where $\proj(x_1 x_2 x_3)(i) = |\{j\in \{1,2,3\} : x_j\leq x_i \}|
$ for $i\in\{1,2,3\}$ and $x_1,x_2,x_3\in [n]$. For instance,
$\proj(127)=\proj(138)=\proj(238)=123$, and
$$(\axbc)\, 491273865 = |\{127,138,238\}| = 3.
$$
Similarly, we also define $(xy\dd z)\,\pi = (z\dd yx)\,\pi^r$,
where $\pi^r$ denotes the reverse of $\pi$, that is, $\pi$ read
backwards.

For any word (finite sequence of letters), $w$, we denote by $|w|$ the
length of $w$, that is, the number of letters in $w$.  A pattern
$\pt=\pt_1\dd\pt_2\dd\cdots\dd\pt_k$ containing exactly $k-1$ dashes
is said to be of type $(|\pt_1|,|\pt_2|,\ldots,|\pt_k|)$.  For
example, the pattern $142\dd 5\dd 367$ is of type $(3,1,3)$, and any
classical pattern of length $k$ is of type
$(\underbrace{1,1,\ldots,1}_{k-1})$.

In \cite{ElNo01} Elizalde and Noy presented the following theorem
regarding the distribution of the number of occurrences of any pattern
of type $(3)$.

\begin{theorem}[Elizalde and Noy \cite{ElNo01}] 
  Let $h(x)=\sqrt{(x-1)(x+3)}$. Then
  \begin{align*}
    \sum_{\pi\in\sym} x^{(\abc)\pi}\frac{t^{|\pi|}}{|\pi|!}
    &=\frac{2h(x)e^{\frac{1}{2}(h(x)-x+1)t}}
    {h(x)+x+1+(h(x)-x-1)e^{h(x)t}},\\
    \sum_{\pi\in\sym} x^{(\bac)\pi}\frac{t^{|\pi|}}{|\pi|!}
    &=\frac{1}{1-\int_0^{t}e^{(x-1)z^2/2}dz}.
  \end{align*}
\end{theorem}

The easy proof of the following proposition can be found in \cite{Cl01}.

\begin{proposition}[Claesson \cite{Cl01}]
  With respect to being equidistributed, the twelve patterns of type
  $(1,2)$ or $(2,1)$ fall into the three classes
  $$
  \begin{array}{l}
  \{\;\axbc,\;\cxba,\;\abxc,\;\cbxa\;\},\\
  \{\;\axcb,\;\cxab,\;\baxc,\;\bcxa\;\},\\
  \{\;\bxac,\;\bxca,\;\acxb,\;\caxb\;\}.
  \end{array}
  $$
\end{proposition}

In the subsequent discussion we refer to the classes of
the proposition above (in the order that they appear) as Class
1, 2 and 3 respectively.

Claesson \cite{Cl01} also gave a solution for the number of permutations
avoiding any pattern of the type $(1,2)$ or $(2,1)$ as follows.

\begin{proposition}[Claesson \cite{Cl01}]\label{claesson1}
  Let $n\in\N$. We have
  $$
  |\sym_n(\pt)| = 
  \begin{cases}
    B_n & \text{ if }\,
    \pt\in\{
    \axbc, \cxba, \abxc, \cbxa,
    \axcb, \cxab, \baxc, \bcxa
    \},\\
    C_n &\text{ if }\,
    \pt\in\{
    \bxac, \bxca, \acxb, \caxb
    \},\\
  \end{cases}
  $$
  where $B_n$ and $C_n$ are the $n$th Bell and Catalan numbers,
  respectively.
\end{proposition}

In particular, since $B_n$ is not $P$-recursive in $n$, this result
implies that for generalized patterns the conjecture that
$s_{\sigma}^r(n)$ is $P$-recursive in $n$ is false for $r=0$ and, 
for example, $\sigma=\axbc$.

This paper is organized as follows. In Section~\ref{class1&2} we find
a recursion for the number of permutations containing exactly one
occurrence of a pattern of Class~1, and we also find an ordinary
generating function for these numbers. We prove these results both
combinatorially and analytically.  Similar results are also obtained
for patterns of Class~2.  In Section~\ref{class3} we give the
distribution of any pattern of Class~3 in the form of a continued
fraction, and we also give explicit formulas for the number of
permutations containing exactly $r$ occurrences of a pattern of
Class~3 when $r\in\{1,2,3\}$.

\section{Counting occurrences of a pattern of Class 1 or 2}\label{class1&2}

\begin{theorem}\label{u1}
  Let $u_1(n)$ be the number of permutations of length $n$ containing
  exactly one occurrence of the pattern $\axbc$ and let $B_n$ be the
  $n$th Bell number. The numbers $u_1(n)$ satisfy the recurrence
  $$ u_1(n+2) = 2 u_1(n+1) + 
  \sum_{k=0}^{n-1}\binom{n}{k}\bigl[u_1(k+1)+B_{k+1}\bigr],
  $$
  whenever $n\geq -1$, with the initial condition $u_1(0)=0$.
\end{theorem}

\begin{proof}
  Each permutation $\pi\in\sym_{n+2}^1(\axbc)$ contains a unique
  subword $abc$ such that $a<b<c$ and $bc$ is a segment of $\pi$. Let
  $x$ be the last letter of $\pi$ and define the sets $\T$, $\T'$, and
  $\T''$ by
  $$ \pi \in
  \begin{cases}
    \T   & \text{ if $x = 2$},\\
    \T'  & \text{ if $x \neq 2$ and $a=1$},\\
    \T'' & \text{ if $x \neq 2$ and $a\neq 1$}.
  \end{cases}
  $$
  Then $\sym_{n+2}^1(\axbc)$ is the disjoint union of $\T$, $\T'$,
  and $\T''$, so 
  $$u_1(n+2)=|\T|+|\T'|+|\T''|.
  $$

  Since removing/adding a trailing $2$ from/to a permutation does
  not affect the number of hits of $\axbc$, we immediately get
  $$|\T| = u_1(n+1).
  $$
  
  For the cardinality of $\T'$ we observe that if $x\neq 2$ and $a=1$
  then $b=2$: If the letter $2$ precedes the letter $1$ then every hit
  of $\axbc$ with $a=1$ would cause an additional hit of $\axbc$ with
  $a=2$ contradicting the uniqueness of the hit of $\axbc$; if $1$
  precedes $2$ then $a=1$ and $b=2$. Thus we can factor any
  permutation $\pi\in \T'$ uniquely in the form $\pi = \sigma 2\tau$,
  where $\sigma$ is $(\axbc)$-avoiding, the letter $1$ is included in
  $\sigma$, and $\tau$ is nonempty and $(\ab)$-avoiding. Owing to
  Proposition~\ref{claesson1} we have showed
  $$|\T'| = \sum_{k=0}^{n-1}\binom{n}{k}B_{k+1}.
  $$
  
  Suppose $\pi\in \T''$. Since $x\neq 2$ and $a\neq 1$ we can factor
  $\pi$ uniquely in the form $\pi = \sigma 1\tau$, where $\sigma$
  contains exactly one occurrence of $\axbc$, the letter $2$ is
  included in $\sigma$, and $\tau$ is nonempty and $(\ab)$-avoiding.
  Consequently,
  $$|\T''| = \sum_{k=0}^{n}\binom{n}{k}u_1(k+1),
  $$
  which completes the proof.
\end{proof}

\begin{example}
  Let us consider all permutations of length $5$ that contain exactly
  one occurrence of $\axbc$, and give a small illustration of the
  proof of Theorem~\ref{U1}. If $\T$, $\T'$ and $\T''$ are defined as
  above then 
  \vspace{1.4ex}
  \newcommand{\ul}{\underline}
  \begin{eqnarray*}
    \T &=& 
      \begin{aligned}[c]
        \mbox{}
        &  \ul{1}\ul{3}\ul{5}4|2\;\;\ul{1}4\ul{3}\ul{5}|2\;\;
        \ul{1}\ul{4}\ul{5}3|2\;\;\ul{1}5\ul{3}\ul{4}|2\;\;
        4\ul{1}\ul{3}\ul{5}|2\;\;5\ul{1}\ul{3}\ul{4}|2\;\;
        \ul{3}\ul{4}\ul{5}1|2
      \end{aligned}\\[1.5ex]
    \T'&=&
      \begin{aligned}[c]
        \mbox{}
        & \ul{1}|\ul{2}\ul{5}43\;\;\ul{1}3|\ul{2}\ul{5}4\;\;
        \ul{1}4|\ul{2}\ul{5}3\;\;\ul{1}43|\ul{2}\ul{5}\;\;
        \ul{1}5|\ul{2}\ul{4}3\;\;\ul{1}53|\ul{2}\ul{4}\\
        & \ul{1}54|\ul{2}\ul{3}\;\;3\ul{1}|\ul{2}\ul{5}4\;\;
        3\ul{1}4|\ul{2}\ul{5}\;\;3\ul{1}5|\ul{2}\ul{4}\;\;
        34\ul{1}|\ul{2}\ul{5}\;\;35\ul{1}|\ul{2}\ul{4}\\
        & 4\ul{1}|\ul{2}\ul{5}3\;\;4\ul{1}3|\ul{2}\ul{5}\;\;
        4\ul{1}5|\ul{2}\ul{3}\;\; 43\ul{1}|\ul{2}\ul{5}\;\;
        45\ul{1}|\ul{2}\ul{3}\;\;5\ul{1}|\ul{2}\ul{4}3\\
        & 5\ul{1}3|\ul{2}\ul{4}\;\;5\ul{1}4|\ul{2}\ul{3}\;\;
        53\ul{1}|\ul{2}\ul{4}\;\;54\ul{1}|\ul{2}\ul{3}
      \end{aligned}\\[1.5ex]
    \T''&=&
      \begin{aligned}[c]
        \mbox{}
        & \ul{2}\ul{3}\ul{4}|15\;\;\ul{2}\ul{3}\ul{5}|14\;\;
        \ul{2}\ul{3}\ul{5}4|1\;\;\ul{2}4\ul{3}\ul{5}|1\;\;
        \ul{2}\ul{4}\ul{5}|13\\
        & \ul{2}\ul{4}\ul{5}3|1\;\;\ul{2}5\ul{3}\ul{4}|1\;\;
        \ul{3}\ul{4}\ul{5}2|1\;\;4\ul{2}\ul{3}\ul{5}|1\;\;
        5\ul{2}\ul{3}\ul{4}|1
      \end{aligned}
  \end{eqnarray*}\\
  where the underlined subword is the unique hit of $\axbc$, and the bar
  indicates how the permutation is factored in the proof of
  Theorem~\ref{U1}.
\end{example}

\begin{theorem}\label{v1}
  Let $v_1(n)$ be the number of permutations of length $n$ containing
  exactly one occurrence of the pattern $\axcb$ and let $B_n$ be the
  $n$th Bell number. The numbers $v_1(n)$ satisfy the recurrence
  $$ v_1(n+1) = v_1(n) + \sum_{k=1}^{n-1}
    \left[
      \binom{n}{k}v_1(k) + \binom{n-1}{k-1}B_k
    \right],
  $$  
  whenever $n\geq 0$, with the initial condition $v_1(0)=0$.
\end{theorem}

\begin{proof}
  Each permutation $\pi\in\sym_{n+2}^1(\axcb)$ contains a unique
  subword $acb$ such that $a<b<c$ and $cb$ is a segment of $\pi$.
  Define the sets $\T$ and $\T'$ by
  $$ \pi \in
  \begin{cases}
    \T   & \text{ if $a = 1$},\\
    \T'  & \text{ if $a \neq 1$}.\\
  \end{cases}
  $$
  Then $\sym_{n+2}^1(\axcb)$ is the disjoint union of $\T$ and $\T'$, so 
  $$v_1(n+2)=|\T|+|\T'|.
  $$
  
  For the cardinality of $\T$ we observe that if $a=1$ then $b=2$: If
  the letter $2$ precedes the letter $1$ or $12$ is a segment of $\pi$
  then every hit of $\axbc$ with $a=1$ would cause an additional hit
  of $\axcb$ with $a=2$ contradicting the uniqueness of the hit of
  $\axbc$; if $1$ precedes $2$ then $a=1$ and $b=2$. Thus we can
  factor $\pi$ uniquely in the form $\pi = \sigma x 2\tau$, where
  $\sigma x$ is $(\axcb)$-avoiding, the letter $1$ is included in
  $\sigma$, and $\tau$ is nonempty and $(\ab)$-avoiding.  Let $\R_n$
  be the set of $(\axcb)$-avoiding permutations of $[n]$ that do not
  end with the letter $1$. Since the letter $1$ cannot be the last
  letter of a hit of $\axcb$, we have, by Proposition~\ref{claesson1},
  that $|\sym_n^0(\axcb)\setminus \R_n|=B_{n-1}$. Consequently, $|\R_n| =
  B_n - B_{n-1}$ and
  \begin{align*}
    |\T| &= \sum_{k=1}^n    \binom{n-1}{k-1} |\R_k| \\
         &= \sum_{k=1}^n    \binom{n-1}{k-1} (B_k-B_{k-1}) \\
         &= \sum_{k=1}^{n-1}\binom{n-1}{k-1} B_k.
  \end{align*} 
  For the last identity we have used the familiar recurrence relation
  $B_{n+1}=\sum_{k=0}^n\binom n k B_k$.

  Suppose $\pi\in \T'$. Since $a\neq 1$ we can factor
  $\pi$ uniquely in the form $\pi = \sigma 1\tau$, where $\sigma$
  contains exactly one occurrence of $\axcb$, and $\tau$ is nonempty
  and $(\ab)$-avoiding. Accordingly,
  $$|\T''| = \sum_{k=0}^{n}\binom{n}{k}v_1(k),
  $$
  which completes the proof.
\end{proof}
            
Let $\pt$ be a pattern of Class 1 or 2. Using combinatorial
reasoning we have found a recursion for the number of permutations
containing exactly one occurrence of the pattern $\pt$
(Theorem~\ref{u1} and \ref{v1}).  More generally, given $r\geq 0$, we
would like to find a recursion for the number of permutations
containing exactly $r$ occurrence of the pattern $\pt$. Using a more
general and analytic approach we will now demonstrate how this (at least
in principle) can be achieved.

Let $S_\pt^r(x)$ be the generating function $S_\pt^r(x)=\sum_n
s_{\pt}^r(n)x^n$.  To find functional relations for $S_\pt^r(x)$ the
following lemma will turn out to be useful.

\begin{lemma}\label{transf}
  If $\{a_n\}$ is a sequence of numbers and $A(x)=\sum_{n\geq 0}
  a_n x^n$ is its ordinary generating function, then, for any $d\geq 0$,
  $$\sum_{n\geq 0} \left[\sum_{j=0}^n \binom{n}{j} a_{j+d}\right]
  x^n= \frac{(1-x)^{d-1}} {x^d} \left[ A\Bigl( \frac x
    {1-x}\Bigr)-\sum_{j=0}^{d-1} a_j \Bigl(\frac x {1-x}\Bigr)^j
  \right].
  $$
\end{lemma}

\begin{proof}
  It is plain that 
  $$
  \sum\limits_{n\geq 0} 
  \Bigl[\sum\limits_{j=0}^n \binom{n}{j} a_{j}\Bigr]
  x^n= \frac 1 {1-x} A\Bigl( \frac x {1-x}\Bigr).
  $$
  See for example \cite[p 192]{GrKnPa94}. On the other hand,   
  $$
  \sum\limits_{n\geq 0} a_{n+d} x^n = 
  \frac 1 {x^d} \Bigl[A(x)-\sum\limits_{j=0}^{d-1} a_j x^j \Bigr].
  $$
  Combining these two identities we get the desired result.
\end{proof}

Define $\sym_n^r(\pt)$ to be the set of permutations $\pi\in S_n$ such
that $(\pt)\pi=r$. Let $s_{\pt}^r(n)= |\sym_n^r(\pt)|$ for $r\geq 0$
and $s_\pt^r(n)=0$ for $r<0$.  Given $b_1,b_2,\ldots,b_k\in\N$, we
also define
$$s_{\pt}^r(n;b_1,b_2,\dots,b_k)=
\#\{a_1 a_2 \cdots a_n\in \sym_n^r(\pt) \mid a_1 a_2 \cdots a_k = b_1
b_2\cdots b_k \}.
$$
As a direct consequence of the above definitions, we have
\begin{equation}\label{eq1}
  s^r_{\pt}(n)=\sum_{j=1}^n s_{\pt}^r(n;j).
\end{equation}

We start by considering patterns that belong to Class 1 and we use
$\abxc$ as a representative of this class. Let us define
\begin{eqnarray*}
u_r(n;b_1,\dots,b_k) &=& s_{\abxc}^r(n;b_1,\dots,b_k),\\
              u_r(n) &=& s_{\abxc}^r(n), \\
              U_r(x) &=& S_{\abxc}^r(x).\\
\end{eqnarray*}

\begin{lemma}\label{tt1}
  Let $n\geq 1$. We have $u_r(n;n-1)=u_r(n;n)=u_r(n-1)$
  and
  $$u_r(n;i)=\sum\limits_{j=1}^{i-1} u_r(n-1;j)
  +\sum\limits_{j=0}^{n-i-1} u_{r-j}(n-1;n-1-j),
  $$ 
  whenever $1 \leq i \leq n-2$.
\end{lemma}

\begin{proof}
  If $a_1 a_2 \cdots a_n$ is any permutation of $[n]$ then  
  $$
  (\abxc)a_1 a_2 \cdots a_n = (\abxc) a_2 a_3\cdots a_n +
  \begin{cases}
    n-a_2  & \text{if } a_1<a_2,\\
    0      & \text{if } a_1>a_2.
  \end{cases}
  $$
  Hence,
  \begin{eqnarray*}
    u_r(n;i) 
    &=& \sum_{j=1}^{i-1} u_r(n;i,j)+\sum_{j=i+1}^n u_r(n;i,j)\\
    &=& \sum_{j=1}^{i-1} u_r(n-1;j)+\sum_{j=i+1}^n u_{r-n+j}(n-1;j-1)\\
    &=& \sum_{j=1}^{i-1} u_r(n-1;j)+\sum_{j=0}^{n-i-1}
    u_{r-j}(n-1;n-1-j).
  \end{eqnarray*}
  For $i=n-1$ or $i=n$ it is easy to see that $u_r(n;i)=u_r(n-1)$.
\end{proof}

Using Lemma \ref{tt1} we quickly generate the numbers $u_r(n)$;
the first few of these numbers are given in Table \ref{tb1}. 
\begin{table}[ht]
  $$
  \begin{array}{rrrrrrrr}
    n\backslash r & 0 & 1 & 2 & 3  & 4  & 5  & 6  \\
    \hline
    0\;\;  & 1    &      &      &      &      &      &      \\
    1\;\;  & 1    &      &      &      &      &      &      \\
    2\;\;  & 2    &      &      &      &      &      &      \\
    3\;\;  & 5    & 1    &      &      &      &      &      \\
    4\;\;  & 15   & 7    &1     &1     &      &      &      \\
    5\;\;  & 52   & 39   &13    &12    &2     &1     &1     \\
    6\;\;  & 203  &211   &112   &103   &41    &24    &17    \\
    7\;\;  & 877  &1168  &843   &811   &492   &337   &238   \\
    8\;\;  & 4140 &6728  &6089  &6273  &4851  &3798  &2956  \\
    9\;\;  & 21147&40561 &43887 &48806 &44291 &38795 &33343 \\
   10\;\;  &115975&256297&321357&386041&394154&379611&355182\\
   \hline
  \end{array}
  $$
  \caption{The number of permutations of length $n$ containing 
    exactly $r$ occurrences of the pattern $12\mbox{-}3$.
    }
  \label{tb1}
\end{table}
Given $r\in\N$ we can also use Lemma \ref{tt1} to find a functional
relation determining $U_r(x)$. Here we present such functional
relations for $r=0,1,2$ and also explicit formulas for $r=0,1$.

Equation $\ref{eq1}$ tells us how to compute $u_r(n)$ if we are given
the numbers $u_r(n;i)$. For the case $r=0$ Lemma~\ref{u0i}, below,
tells us how to do the converse.

\begin{lemma}\label{u0i}
  If $1 \leq i \leq n-2$ then
  $$u_0(n;i)=\sum_{j=0}^{i-1} \binom{i-1}{j} u_0(n-2-j).$$
\end{lemma}

\begin{proof}
  For $n=1$ the identity is trivially true. Assume the identity is
  true for $n=m$. We have
  \begin{align*}
    u_0(m+1;i) 
    &= \sum_{j=1}^{i-1}u_0(m;j) + u_0(m-1) 
    && \text{by Lemma~\ref{tt1}}\\
    &= \sum_{j=1}^{i-1}\sum_{k=0}^{j-1}\binom{j-1}{k}u_0(m-2-k)+u_0(m-1)
    && \parbox{17ex}{by the induction hypothesis}\\
    &= \sum_{j=1}^{i-1}\sum_{k=j-1}^{i-2}\binom{k}{j-1}u_0(m-1-j).
  \end{align*}
  Using the familiar equality $\binom 1 k + \binom 2 k + \cdots +
  \binom n k = \binom{n+1}{k+1}$ we then get
  $$u_0(m+1;i) = \sum_{j=1}^{i-1}\binom{i-1}{j} u_0(m-1-j).
  $$
  Thus the identity is true for $n=m+1$ and by the principle of
  induction the desired identity is true for all $n\geq 1$.
\end{proof}

The following proposition is a direct consequence of
Proposition~\ref{claesson1}. However, we give a different proof. The
proof is intended to illustrate the general approach. It is advisable
to read this proof before reading the proof of
Theorem~\ref{u1}$^{\prime}$ below.

\begin{proposition}\label{U0}
  The ordinary generating function for the number of
  $(\abxc)$-avoiding permutations of length $n$ is
  $$U_0(x)=\sum_{k\geq 0} \frac{x^k}{(1-x)(1-2x)\cdots(1-kx)}.
  $$
\end{proposition}

\begin{proof}
  We have
  \begin{align*}
    u_0(n) 
    &= \sum_{k=1}^n u_0(n;k) 
    && \text{by Equation \ref{eq1}}\\
    &= 2u_0(n-1)+\sum_{i=1}^{n-2}\sum_{j=0}^{i-1}\binom{i-1}{j}u_0(n-2-j)
    && \text{by Lemma~\ref{tt1} and \ref{u0i}}\\
    &= u_0(n-1)+\sum_{i=0}^{n-2}\binom{n-2}{i}u_0(n-1-i)
    && \text{by \scriptsize{
        $\sum\limits_{i=k}^n\binom i k=\binom{n+1}{k+1}$ 
        }}\\
    &= u_0(n-1)+\sum_{i=0}^{n-2}\binom{n-2}{i}u_0(i+1).
  \end{align*}
  Therefore, by Lemma \ref{transf}, we have
  $$U_0(x)=xU_0(x)+1-x+xU_0\left( \frac{x}{1-x} \right)\!,
  $$
  which is equivalent to 
  $$U_0(x)=1+\frac{x}{1-x}U_0\left( \frac{x}{1-x} \right)\!.
  $$
  An infinite number of applications of this identity concludes the proof.
\end{proof}

We now derive a formula for $U_1(x)$ that is somewhat similar to the
one for $U_0(x)$. The following lemma is a first step in this
direction.

\begin{lemma}\label{u1i}
  If $1 \leq i \leq n-2$ then
  $$u_1(n;i)=\sum_{j=0}^{i-1} \binom{i-1}{j} u_1(n-2-j)+u_0(n;i).$$
\end{lemma}

\begin{proof}
  For $n=1$ the identity is trivially true. Assume the identity is
  true for $n=m$. Lemma \ref{tt1} and the induction hypothesis imply
  \begin{eqnarray*}
    u_1(m+1;i) 
    &=& \sum_{j=1}^{i-1}u_1(m;j) + u_1(m-1) + u_0(m-1) \\
    &=& \sum_{j=0}^{i-1}\binom{j-1}{k} u_1(m-1-j) + 
    \sum_{j=1}^{i-1}u_0(m;j)+u_0(m-1).
  \end{eqnarray*}
  In addition, Lemma~\ref{u0i} implies
  \begin{eqnarray*}
   u_0(m+1;i) 
   &=& \sum_{j=1}^{i-1}\sum_{k=0}^{j-1}\binom{j-1}{k}u_0(n-2-k) + u_0(n-1)\\
   &=& \sum_{j=0}^{i-1}\binom{i-1}{j}u_0(n-1-j)\\ 
   &=& \sum_{j=1}^{i-1}u_0(m;j)+u_0(m-1).
  \end{eqnarray*}
  Thus the identity is true for $n=m+1$ and by the principle of induction
  the desired identity is true for all $n\geq 1$.
\end{proof}

Next, we rediscover Theorem~\ref{u1}.

\begin{thmu1}
  Let $u_1(n)$ be the number of permutations of length $n$ containing
  exactly one occurrence of the pattern $\abxc$ and let $B_n$ be the
  $n$th Bell number. The numbers $u_1(n)$ satisfy the recurrence
  $$ u_1(n+2) = 2 u_1(n+1) + 
  \sum_{k=0}^{n-1}\binom{n}{k}\bigl[u_1(k+1)+B_{k+1}\bigr],
  $$
  whenever $n\geq -1$, with the initial condition $u_1(0)=0$.
\end{thmu1}

\begin{proof}
  Similarly to the proof of Proposition~\ref{U0}, we use
  Equation~\ref{eq1}, Lemma~\ref{tt1}, \ref{u0i}, and \ref{u1i} to get
  \begin{eqnarray*}
    u_1(n) 
    &=& 2 u_1(n-1) + 
    \sum_{i=1}^{n-2}\left[
      \sum_{j=0}^{i-1}\binom{i-1}{j}u_1(n-2-j)+u_0(n;i)
    \right]\\
    &=& 2 u_1(n-1) + 
    \sum_{i=1}^{n-2}\sum_{j=0}^{i-1}\binom{i-1}{j}\bigl(
      u_1(n-2-j)+u_0(n-2-j)
    \bigr)\\
    &=& u_1(n-1) - u_0(n-1) + 
    \sum_{i=0}^{n-2}\binom{n-2}{i}\bigl(u_1(i+1)+u_0(i+1)\bigr)\\
    &=& 2u_1(n-1) + 
    \sum_{i=0}^{n-3}\binom{n-2}{i}\bigl(u_1(i+1)+u_0(i+1)\bigr).
  \end{eqnarray*}
\end{proof}

\begin{corollary}\label{U1}
  The ordinary generating function, $U_1(x)$, for the number of
  permutations of length $n$ containing exactly one occurrence of the
  pattern $\abxc$ satisfies the functional equation
  $$U_1(x)=\frac{x}{1-x}\biggl(
  U_1\Bigl(\frac{x}{1-x}\Bigr)+
  U_0\Bigl(\frac{x}{1-x}\Bigr)-U_0(x)
  \biggr).
  $$
\end{corollary}

\begin{proof} 
  The result follows from Theorem~\ref{u1} together with
  Lemma~\ref{transf}.
\end{proof}

\begin{corollary}
  The ordinary generating function for the number of permutations of
  length $n$ containing exactly one occurrence of the pattern $\abxc$ is
  $$
  U_1(x) = \sum_{n\geq 1}\frac x {1-nx} 
  \sum_{k\geq 0} \frac {kx^{k+n}} {(1-x)(1-2x)\cdots(1-(k+n)x)}.\\
  $$
\end{corollary}

\begin{proof}
  We simply apply Corollary~\ref{U1} an infinite number of times and in
  each step we perform some rather tedious algebraic manipulations.
\end{proof}

\begin{theorem}\label{U2}
  The ordinary generating function, $U_2(x)$, for the number of
  permutations of length $n$ containing exactly two occurrences of the
  pattern $\abxc$ satisfies the functional equation
  $$
  U_2(x)=\frac{x}{(1-x)^2(1-2x)}\biggl(\;
    \begin{aligned}[t]
      &U_2\Bigl(\frac{x}{1-x}\Bigr)-(1-x) U_2(x)+\\
      &U_1\Bigl(\frac{x}{1-x}\Bigr)-(1-x)^2 U_1(x)+\\
      &U_0\Bigl(\frac{x}{1-x}\Bigr)-(1-x)^2 U_0(x)
      \;\;\biggr).
    \end{aligned}
  $$
\end{theorem}

\begin{proof} 
  The proof is similar to the proofs of Lemma~\ref{u1i},
  Theorem~\ref{u1}' and Corollary~\ref{U1}, and we only sketch it here.

 Lemma \ref{tt1} yields
  \begin{align*}
    u_2(n;n)&=u_2(n-1)\\
    u_2(n;n-1)&=u_2(n-1)\\
    u_2(n;n-2)&=u_2(n-1)-u_2(n-2)+u_1(n-2)
  \end{align*}
  and, by means of induction,
  $$u_2(n;i)=u_1(n;i)+u_0(n;i)-u_0(n-1;i)+
  \sum_{j=0}^{i-1}\binom{i-1}{j}u_2(n-2-j),
  $$
  whenever $1\leq i\leq n-3$. Therefore, $u_2(0)=u_2(1)=u_2(2)=0$
  and
  \begin{multline*}
    u_2(n) =3u_2(n-1)-u_2(n-2)+u_1(n-2)+\\\sum_{i=1}^{n-3}
    \binom{n-3}{i}(u_2(n-1-i)+u_1(n-1-i)+u_0(n-1-i)-u_0(n-2-i)).
  \end{multline*}
  whenever $n\geq 3$. Thus, the result follows from Lemma~\ref{transf}.
\end{proof}


We now turn our attention to patterns that belong to Class~2 and we
use $\bcxa$ as a representative of this class.  The results found
below regarding the $\bcxa$ pattern are very similar to the ones
previously found for the $\abxc$ pattern, and so are the proofs;
therefore we choose to omit most of the proofs. However, we give the
necessary lemmas from which the reader may construct her/his own proofs.

Define
\begin{eqnarray*}
v_r(n;b_1,\dots,b_k) &=& s_{\bcxa}^r(n;b_1,\dots,b_k),\\
              v_r(n) &=& s_{\bcxa}^r(n), \\
              V_r(x) &=& S_{\bcxa}^r(x).\\
\end{eqnarray*}
If $a_1 a_2 \cdots a_n$ is any permutation of $[n]$ then  
$$
(\bcxa)a_1 a_2 \cdots a_n = (\bcxa) a_2 a_3\cdots a_n +
\begin{cases}
  a_1-1  & \text{if } a_1<a_2,\\
  0      & \text{if } a_1>a_2.
\end{cases}
$$

\begin{lemma}\label{tt2}
  Let $n\geq 1$. We have $v_r(n;1)=v_r(n;n)=v_r(n-1)$
  and
  $$v_r(n;i)=\sum_{j=1}^{i-1} v_r(n-1;j)+\sum_{j=i}^{n-1}v_{r-i+1}(n-1;j),
  $$ 
  whenever $2 \leq i \leq n-1$.
\end{lemma}

Using Lemma \ref{tt2} we quickly generate the numbers $v_r(n)$;
the first few of these numbers are given in Table \ref{tb2}. 

\begin{table}[ht]
  $$
  \begin{array}{rrrrrrrr}
    n\backslash r & 0  & 1  & 2  & 3  & 4  & 5  & 6  \\
    \hline
    0\;\; & 1    &      &      &      &      &      &      \\
    1\;\; & 1    &      &      &      &      &      &      \\
    2\;\; & 2    &      &      &      &      &      &      \\
    3\;\; & 5    & 1    &      &      &      &      &      \\
    4\;\; & 15   & 6    &3     &      &      &      &      \\
    5\;\; & 52   & 32   &23    &10    &3     &      &      \\
    6\;\; & 203  &171   &152   &98    &62    &22    &11    \\
    7\;\; & 877  &944   &984   &791   &624   &392   &240   \\
    8\;\; & 4140 &5444  &6460  &6082  &5513  &4302  &3328  \\
    9\;\; & 21147&32919 &43626 &46508 &46880 &41979 &36774 \\
   10\;\; & 115975&208816&304939&360376&396545&393476&377610\\
   \hline
  \end{array}
  $$
  \caption{The number of permutations of length $n$ containing 
    exactly $r$ occurrences of the pattern $23\mbox{-}1$.
    }
  \label{tb2}
\end{table}

\begin{lemma}\label{v0i}
  If $2 \leq i \leq n-1$ then
  $$v_0(n;i)=\sum_{j=0}^{i-2} \binom{i-2}{j} v_0(n-2-j).$$
\end{lemma}

\begin{proposition}\label{V0}
  The ordinary generating function for the number of
  $(\bcxa)$-avoiding permutations of length $n$ is
  $$V_0(x)=\sum_{k\geq 0} \frac{x^k}{(1-x)(1-2x)\cdots(1-kx)}.
  $$
\end{proposition}

\begin{lemma}\label{v1i}
  If $2 \leq i \leq n-1$ then
  $$v_1(n;i)=\sum_{j=0}^{i-2} \binom{i-2}{j}
  v_1(n-2-j)+v_0(n;i-1)-v_0(n-1,i-1).
  $$
\end{lemma}

\begin{thmv1}
  Let $v_1(n)$ be the number of permutations of length $n$ containing
  exactly one occurrence of the pattern $\bcxa$ and let $B_n$ be the
  $n$th Bell number. The numbers $v_1(n)$ satisfy the recurrence
  $$ v_1(n+1) = v_1(n) + \sum_{k=1}^{n-1}
    \left[
      \binom{n}{k}v_1(k) + \binom{n-1}{k-1}B_k
    \right],
  $$  
  whenever $n\geq 0$, with the initial condition $v_1(0)=0$.
\end{thmv1}

\begin{corollary}\label{V1}
  The ordinary generating function for the number of
  permutations of length $n$ containing exactly one occurrence of the
  pattern $\bcxa$ satisfies the functional equation
  $$V_1(x)=\frac{x}{1-x}V_1\Bigl(\frac{x}{1-x}\Bigr)+ 
  x\biggl(
    V_0\Bigl(\frac{x}{1-x}\Bigr)-V_0(x)
  \biggr)\!.
  $$
\end{corollary}

\begin{corollary}
  The ordinary generating function for the number of permutations of
  length $n$ containing exactly one occurrence of the pattern $\bcxa$ is
  $$
  V_1(x) = \sum_{n\geq 1}\frac x {1-(n-1)x} 
  \sum_{k\geq 0} \frac {kx^{k+n}} {(1-x)(1-2x)\cdots(1-(k+n)x)}.\\
  $$
\end{corollary}

\begin{theorem}
  The ordinary generating function, $V_2(x)$, for the number of
  permutations of length $n$ containing exactly two occurrences of the
  pattern $\bcxa$ satisfies the functional equation
  $$
  V_2(x)=\frac x{1-x}\biggl(V_2\Bigl(\frac
  x{1-x}\Bigr)+(1-2x)V_1\Bigl(\frac
  x{1-x}\Bigr)+(1-3x+x^2)V_0\Bigl(\frac x{1-x}\Bigr)\biggr)-x+x^2
  $$
\end{theorem}

\begin{proof}
  By Lemma 5
  \begin{align*}
    v_2(n;n)&=v_2(n-1)\\
    v_2(n;1)&=v_2(n-1)\\
    v_2(n;2)&=v_2(n-2)+v_1(n-1)-v_1(n-2)\\
    v_2(n;3)&=v_2(n-2)+v_2(n-3)+v_1(n-2)-v_1(n-3)+\\
    &\hspace{30ex}+v_0(n-1)-v_0(n-2)-v_0(n-3)
  \end{align*}
  and, by means of induction,
  $$
  v_2(n;i)=\sum_{j=0}^{i-2}\binom{i-2}{j}v_2(n-2-j)+
  v_1(n;i-1)+v_1(n-1;i-1)-v_0(n-1;i-2)
  $$
  for $n-1\geq i\geq 4$. Thus $v_2(0)=v_2(1)=v_2(2)=0$ and for all $n\geq 3$
  \begin{align*}
    v_2(n)=v_2(n-1)+&\sum_{j=0}^{n-2} \binom{n-2}{j} v_2(n-1-j)+\\
    +&\sum_{j=0}^{n-3}\binom{n-3}{j}\bigl(v_1(n-1-j)-v_1(n-2-j)\bigr)+\\
    +&\sum_{j=0}^{n-4}\binom{n-4}{j}\bigl(v_0(n-1-j)-v_0(n-2-j)-v_0(n-3-j)\bigr).
  \end{align*}
  The result now follows from Lemma~\ref{transf}.
\end{proof}

\section{Counting occurrences of a pattern of Class 3}\label{class3}

We choose $\bxac$ as our representative for Class~3 and we define
$w_r(n)$ as the number of permutations of length $n$ containing exactly
$r$ occurrences of the pattern $\bxac$.  We could apply the analytic
approach from the previous section to the problem of determining
$w_r(n)$.  However, a result by Clarke, Steingr{\'\i}msson and Zeng
\cite[Corollary~11]{ClStZe97} provides us with a better option.

\begin{theorem}\label{CSZ}
  The following Stieltjes continued fraction expansion holds
  $$ 
  \sum_{\pi\in\sym} x^{1+(\ab)\pi } y^{(\ba)\pi} 
  p^{(\bxca)\pi} q^{(\caxb)\pi} t^{|\pi|} = 
  \cfrac{1}{1 -
    \cfrac{x \pq 1 t}{1 - 
      \cfrac{y \pq 1 t}{1 -
        \cfrac{x \pq 2 t}{1 - 
          \cfrac{y \pq 2 t}{\quad\ddots
            }}}}}
  $$
  where $\pq n = q^{n-1}+pq^{n-2}+\cdots+p^{n-2}q+p^{n-1}$.
\end{theorem}

\begin{proof}
  In \cite[Corollary~11]{ClStZe97} Clarke, Steingr{\'\i}msson and Zeng
  derived the following continued fraction expansion
  $$ 
  \sum_{\pi\in\sym}  
  y^{\des \pi} p^{\res\pi} q^{\ddif\pi} t^{|\pi|} = 
  \cfrac{1}{1 -
    \cfrac{         \p 1 t}{1 - 
      \cfrac{y q     \p 1 t}{1 -
        \cfrac{q      \p 2 t}{1 - 
          \cfrac{y q^2 \p 2 t}{\quad\ddots
            }}}}}
  $$
  where $\p n = 1+p+\cdots+p^{n-1}$. We refer the reader to
  \cite{ClStZe97} for the definitions of $\ddif$ and $\res$. However,
  given these definitions, it is easy to see that $\res = (\bxca)$ and
  $\ddif = (\ba)+(\bxca)+(\caxb)$. Moreover, $\des = (\ba)$ and
  $|\pi|=1+(\ab)\pi+(\ba)\pi$. Thus, substituting $y (xq)^{-1}$ for
  $y$, $p q^{-1}$ for $p$, and $xt$ for $t$, we get the desired result.
\end{proof}

The following corollary is an immediate consequence of Theorem~\ref{CSZ}.

\begin{corollary}\label{cf}
  The bivariate ordinary generating function for the distribution of
  occurrences of the pattern $\bxac$ admits the Stieltjes continued
  fraction expansion
  $$
  \sum_{\pi\in\sym} p^{(\bxac)\pi} t^{|\pi|}=   
  \cfrac{1}{1 -
    \cfrac{ \p 1 t}{1 - 
      \cfrac{\p 1 t}{1 -
        \cfrac{\p 2 t}{1 - 
          \cfrac{\p 2 t}{\quad\ddots
            }}}}}
  $$
  where $\p n = 1+p+\cdots+p^{n-1}$ 
\end{corollary}

Using Corollary~\ref{cf} we quickly generate the numbers $w_r(n)$;
the first few of these numbers are given in Table \ref{tb3}. 

\begin{table}[ht]
  $$
  \begin{array}{rrrrrrrr}
    n\backslash r & 0& 1 & 2  & 3  & 4  & 5  & 6 \\
    \hline
    0\;\;  & 1    &      &      &      &      &      &      \\
    1\;\;  & 1    &      &      &      &      &      &      \\
    2\;\;  & 2    &      &      &      &      &      &      \\
    3\;\;  & 5    & 1    &      &      &      &      &      \\
    4\;\;  & 14   & 8    &2     &      &      &      &      \\
    5\;\;  & 42   & 45   &25    &7     &1     &      &      \\
    6\;\;  & 132  &220   &198   &112   &44    &12    &2     \\
    7\;\;  & 429  &1001  &1274  &1092  &700   &352   &140   \\
    8\;\;  & 1430 &4368  &7280  &8400  &7460  &5392  &3262  \\
    9\;\;  & 4862 &18564 &38556 &56100 &63648 &59670 &47802 \\
    10\;\; &16796 &77520 &193800&341088&470934&541044&535990\\
    \hline
  \end{array}
  $$
  \caption{The number of permutations of length $n$ containing 
    exactly $r$ occurrences of the pattern $2\mbox{-}13$.
    }
  \label{tb3}
\end{table}

\begin{corollary}
  The number of $(\bxac)$-avoiding permutations of length $n$ is
  $$w_0(n) = \frac{1}{n+1}\binom{2n}{n}.
  $$
\end{corollary}

\begin{proof}
  This result is explicitly stated in Proposition~\ref{claesson1},
  but it also follows from Corollary~\ref{cf} by putting $p=0$.
\end{proof}

\begin{corollary}\label{w1}
  The number of permutations of length $n$ containing
  exactly one occurrence of the pattern $\bxac$ is
  $$w_1(n) = \binom{2n}{n-3}.
  $$
\end{corollary}

\begin{proof}
  For $m>0$ let
  $$W(p,t;m) = 
  \cfrac{1}{1 -
    \cfrac{ \p m t}{1 - 
      \cfrac{\p m t}{1 -
        \cfrac{\p {m+1} t}{1 - 
          \cfrac{\p {m+1} t}{\quad\ddots
            }}}}}
  $$
  Note that
  $$W(p,t;m) = 
  \cfrac{1}{1 -
    \cfrac{ \p m t}{1 - \p m t W(p,t;m+1).
      }} 
  $$
  Assume $m>1$. Differentiating $W(p,t;m)$ with respect to $p$ and
  evaluating the result at $p=0$ we get
  $$D_pW(p,t;m)\big{|}_{p=0}=tC(t)^3+t^2C(t)^5+t^2C(t)^4
  D_pW(p,t;m+1)\big{|}_{p=0}
  $$
  where $C(t)=W(0,t,1)$ is the generating function for the Catalan
  numbers. Applying this identity an infinite number of times we get
  $$D_pW(p,t,m)\big{|}_{p=0}
  =tC(t)^3+t^2C(t)^5+t^3C(t)^{7}+\cdots=\frac{t C(t)^3}{1-tC(t)^2}.
  $$
  On the other hand, $D_pW(p,t;1)\big{|}_{p=0}=t^2C(t)^4 D_p
  W(p,t;2)\big{|}_{p=0}$. Combining these two identities we get
  $$D_pW(p,t;1)\big{|}_{p=0} = \frac{t^3C(t)^7}{1-tC(t)^2}.
  $$
  Since $\sum_{n\geq 0}w_1(n) t^n = D_pW(p,t;1)\big{|}_{p=0}$ the
  proof is completed on extracting coefficients in the last identity.
\end{proof}

The proofs of the following two corollaries are similar to the proof
of Corollary~\ref{w1} and are omitted.

\begin{corollary}
  The number of permutations of length $n$ containing exactly two
  occurrences of the pattern $\bxac$ is
  $$ w_2(n) = \frac{n(n-3)}{2(n+4)}\binom{2n}{n-3}.
  $$
\end{corollary}

\begin{corollary}
  The number of permutations of length $n$ containing
  exactly three occurrences of the pattern $\bxac$ is  
  $$ w_3(n) = \frac{1}{3}\binom{n+2}{2}\binom{2n}{n-5}.
  $$
\end{corollary}

As a concluding remark we note that there are many questions left to
answer. What is, for example, the formula for $w_k(n)$ in general?
What is the combinatorial explanation of $n
s^1_{\axbxc}(n)=3s^1_{\bxac}(n)$ and
$$(n+3)(n+2)(n+1)s^1_{\bxac}(n)=2n(2n-1)(2n-2)s^1_{\bxaxc}(n)?
$$
In addition, Corollary~\ref{w1} obviously is in need of a
combinatorial proof.

\bibliographystyle{plain}


\begin{thebibliography}{10}

\bibitem{AlFr00}
N.~Alon and E.~Friedgut.
\newblock On the number of permutations avoiding a given pattern.
\newblock {\em J. Combin. Theory Ser. A}, 89(1):133--140, 2000.

\bibitem{At99}
M.~D. Atkinson.
\newblock Restricted permutations.
\newblock {\em Discrete Math.}, 195(1-3):27--38, 1999.

\bibitem{BaSt00}
E.~Babson and E.~Steingr\'{\i}msson.
\newblock Generalized permutation patterns and a classification of the
  {M}ahonian statistics.
\newblock {\em S\'em. Lothar. Combin.}, 44:Art. B44b, 18 pp. (electronic),
  2000.

\bibitem{Bo97b}
M.~B{\'o}na.
\newblock Exact enumeration of $1342$-avoiding permutations: a close link with
  labeled trees and planar maps.
\newblock {\em J. Combin. Theory Ser. A}, 80(2):257--272, 1997.

\bibitem{Bo97c}
M.~B{\'o}na.
\newblock The number of permutations with exactly $r$ $132$-subsequences is
  ${P}$-recursive in the size!
\newblock {\em Adv. in Appl. Math.}, 18(4):510--522, 1997.

\bibitem{Bo97a}
M.~B{\'o}na.
\newblock Permutations avoiding certain patterns: the case of length $4$ and
  some generalizations.
\newblock {\em Discrete Math.}, 175(1-3):55--67, 1997.

\bibitem{Bo98}
M.~B{\'o}na.
\newblock Permutations with one or two $132$-subsequences.
\newblock {\em Discrete Math.}, 181(1-3):267--274, 1998.

\bibitem{ChWe99}
T.~Chow and J.~West.
\newblock Forbidden subsequences and {C}hebyshev polynomials.
\newblock {\em Discrete Math.}, 204(1-3):119--128, 1999.

\bibitem{Cl01}
A.~Claesson.
\newblock Generalized pattern avoidance.
\newblock {\em European J. Combin.}, 22(7):961--971, 2001.

\bibitem{ClStZe97}
R.J. Clarke, E.~Steingr{\'\i}msson, and J.~Zeng.
\newblock New {E}uler-{M}ahonian statistics on permutations and words.
\newblock {\em Adv. in Appl. Math.}, 18(3):237--270, 1997.

\bibitem{ElNo01}
S.~Elizalde and M.~Noy.
\newblock Enumeration of subwords in permutations.
\newblock In {\em Formal power series and algebraic combinatorics (Tempe,
  2001)}, pages 179--189. Arizona State University, 2001.

\bibitem{GrKnPa94}
R.~L. Graham, D.~E. Knuth, and O.~Patashnik.
\newblock {\em Concrete mathematics}.
\newblock Addison-Wesley Publishing Company, Reading, MA, second edition, 1994.

\bibitem{Ma00}
T.~Mansour.
\newblock Permutations containing and avoiding certain patterns.
\newblock In {\em Formal power series and algebraic combinatorics (Moscow,
  2000)}, pages 704--708. Springer, Berlin, 2000.

\bibitem{MaVa01}
T.~Mansour and A.~Vainshtein.
\newblock Counting occurrences of $132$ in a permutation.
\newblock {\em To appear in: Adv. Appl. Math.}, 2001.

\bibitem{No96}
J.~Noonan.
\newblock The number of permutations containing exactly one increasing
  subsequence of length three.
\newblock {\em Discrete Math.}, 152(1-3):307--313, 1996.

\bibitem{NoZe96}
J.~Noonan and D.~Zeilberger.
\newblock The enumeration of permutations with a prescribed number of
  ``forbidden'' patterns.
\newblock {\em Adv. in Appl. Math.}, 17(4):381--407, 1996.

\bibitem{Ro99}
A.~Robertson.
\newblock Permutations containing and avoiding $123$ and $132$ patterns.
\newblock {\em Discrete Math. Theor. Comput. Sci.}, 3(4):151--154 (electronic),
  1999.

\bibitem{SiSc85}
R.~Simion and F.~W. Schmidt.
\newblock Restricted permutations.
\newblock {\em European J. Combin.}, 6(4):383--406, 1985.

\bibitem{St94}
Z.~Stankova.
\newblock Forbidden subsequences.
\newblock {\em Discrete Math.}, 132(1-3):291--316, 1994.

\bibitem{St96}
Z.~Stankova.
\newblock Classification of forbidden subsequences of length $4$.
\newblock {\em European J. Combin.}, 17(5):501--517, 1996.

\bibitem{We95}
J.~West.
\newblock Generating trees and the {C}atalan and {S}chr\"oder numbers.
\newblock {\em Discrete Math.}, 146(1-3):247--262, 1995.

\end{thebibliography}

\end{document}